\documentclass[12pt]{amsart}
\usepackage{amsmath,amssymb,amsthm,amscd}
\newtheorem{formula}{}[section]
\newtheorem{proposition}[formula]{Proposition}
\newtheorem{corollary}[formula]{Corollary}
\newtheorem{lemma}[formula]{Lemma}
\newtheorem{theorem}[formula]{Theorem}
\theoremstyle{definition}
\newtheorem{definition}[formula]{Definition}
\newtheorem{example}[formula]{Example}
\theoremstyle{remark}
\newtheorem*{remark}{Remark}

\def\A{{\mathcal A}}

\begin{document}

\title[Massey products in Lie algebra cohomology]
{Massey products in graded Lie algebra cohomology}
\subjclass{55S30, 17B56, 17B70, 17B10}
\author{Dmitri Millionschikov}
\address{Moscow State University, Department of Mathematics and Mechanics,
Leninskie gory, 119899 Moscow, Russia}
\email{million@mech.math.msu.su}
\date{December 25, 2005}
\keywords{Massey products, graded Lie algebras, formal connection,
Maurer-Cartan equation, representation, cohomology}
\thanks{The research of the author was partially supported by grants
RFBR 05-01-01032 and "Russian Scientific Schools"}
\begin{abstract}
We discuss Massey products in
a ${\mathbb N}$-graded Lie algebra cohomology.
One of the main examples is so-called "positive part"
$L_1$ of the Witt algebra $W$. Buchstaber conjectured that
$H^*(L_1)$ is generated with respect to non-trivial Massey
products by $H^1(L_1)$. Feigin, Fuchs and Retakh represented
$H^*(L_1)$ by trivial Massey products and the second part of
the Buchstaber conjecture is still open.
We consider an associated graded algebra $\mathfrak{m}_0$
of $L_1$ with respect to the filtration by
the its descending central series and prove
that $H^*(\mathfrak{m}_0)$ is generated with respect to non-trivial Massey
products by $H^1(\mathfrak{m}_0)$.
\end{abstract}
\date{}

\maketitle

\section*{Introduction}
In the last thirty years Massey products in cohomology
have found a lot of interesting
applications in topology and geometry.
The existence of non-trivial Massey products in
$H^*(M,{\mathbb R})$ is an obstruction for a manifold $M$ to be K\"ahler
K\"ahler manifolds are
formal \cite{DGMS}, i.e. their real homotopy types are
completely determined by their real cohomology algebras. In their
turn formal spaces have trivial Massey products.

An important feature of Massey products is the following: a
blow-up of a symplectic manifold $M$ along its submanifold $N$
inherits non-trivial Massey products from $M$ to $N$ \cite{BaTa}.
This idea was used by McDuff \cite{McD} in her construction of
simply connected symplectic manifold with no K\"ahler structures.
The Massey products and examples of symplectic manifolds with no
K\"ahler structures was the subject in \cite{CFG}. Babenko and
Taimanov considered an interesting family of symplectic
nilmanifolds related to graded finite-dimensional quotients of a
"positive part" $L_1$ of the Witt algebra $W$. Applying the
symplectic blow-up procedure they constructed examples of simply
connected non-formal symplectic manifolds in dimensions $\ge 10$
\cite{BaTa0}. In fact it is possible to use another graded Lie
algebra $\mathfrak{m}_0$ instead of $L_1$ \cite{Mill}. A few time
ago an $8$-dimensional example was constructed by Fernandez and
Mu\~noz  \cite{FeMu} using another techniques.

An initial data of almost all examples of non-formal symplectic
manifolds are nilmanifolds related
to positively graded Lie algebras and some non-trivial triple
Massey products in their cohomology. The present article is devoted
to the study of $n$-fold calssical Massey products in the
cohomology of ${\mathbb N}$-graded Lie algebras. Although
we started our introduction with finite-dimensional examples
related to some nilmanifolds we will focuss our attention to
two infinite dimensional ${\mathbb N}$-graded Lie algebras
$L_1$ and its associated graded (with respect to the filtration
by the ideals $C^kL_1$ of the descending central series) Lie algebra
${\rm gr}_C L_1 \cong \mathfrak{m}_0$.
The reason of this interest comes from the relation
of the algebra $L_1$ to the Landweber-Novikov algebra
in the complex cobordisms theory that was discovered by Buchstaber and
Shokurov in $70$-th \cite{BuSho}. In that time the algebra
$L_1$ attracted a lot of researchers \cite{Fu} and the computation
of its cohomology by Goncharova \cite{G} is one of the most technicaly
complicated results in homology algebra.

The cohomology algebra
$H^*(L_1)$ has a trivial multiplication.
Buchstaber conjectured that the algebra
$H^*(L_1)$ is generated with respect to the Massey products
by $H^1(L_1)$, moreover all corresponding Massey
products can be chosen non-trivial.
The first part of Buchstaber's conjecture was proved
by Feigin, Fuchs and Retakh \cite{FeFuRe}. But they represented
the cohomology classes from $H^*(L_1)$ by means of trivial Massey products.
Twelve years later Artel'nykh represented some
of Goncharova's cocycles by non-trivial Massey products,
but his brief article does not contain the proofs.
Thus one may conclude that the original Buchstaber conjecture
is still open.

We recall the necessary information
on the cohomology of graded Lie algebras
in the first Section and study two main examples
$H^*(L_1)$ \cite{G} and $H^*(\mathfrak{m}_0)$ \cite{FialMill}
in the Section \ref{L_1_m_0}. In the Section \ref{Massey_S} we present
May's approach to the definition of Massey products, his notion of formal
connection developped by Babenko and Taimanov in \cite{BaTa} for
Lie algebras, we introduce also the notion of equivalent Massey products.
The analogy with the classical Maurer-Cartan equation is
especially transparent in the case of Massey products
of $1$-dimensional cohomology classes
$\langle \omega_1,\dots,\omega_n\rangle$. The relation
of this special case to the representations theory
was discovered in \cite{FeFuRe}, \cite{Dw}.
We discuss it in the Section \ref{g-modules}.
Following
\cite{FeFuRe} we consider Massey products
$\langle \omega_1,\dots,\omega_n,\Omega\rangle$, where
$\omega_1,\dots,\omega_n$ are closed $1$-forms and
$\Omega$ is a closed $q$-form. These Massey  products are
related to the differentials of Feigin-Fuchs spectral sequence
\cite{FeFuRe}.

The main result of the present article is the Theorem \ref{main_th}
stating that
{\it the cohomology algebra  $H^*(\mathfrak{m}_0)$  is generated
with respect to the non-trivial Massey products
by $H^1(\mathfrak{m}_0)$}. Another important
result is the Theorem \ref{1_Massey} that contains a list of
equivalency classes of trival Massey products of $1$-cohomology classes.
We show that it is related to
Benoist's classification \cite{Ben} of indecomposable
finite-dimensional thread modules over the Lie algebra $\mathfrak{m}_0$.

\section{Cohomology of ${\mathbb N}$-graded Lie algebras}

Let $\mathfrak{g}$ be a Lie algebra over ${\mathbb K}$ and
$\rho: \mathfrak{g} \to \mathfrak{gl}(V)$ its linear representation
(or in other words $V$ is a $\mathfrak{g}$-module).
We denote by $C^q(\mathfrak{g},V)$
the space of $q$-linear skew-symmetric mappings of $\mathfrak{g}$ into
$V$. Then one can consider an algebraic complex:

$$
\begin{CD}
V @>{d_0}>>
C^1(\mathfrak{g}, V) @>{d_1}>> C^2(\mathfrak{g}, V) @>{d_2}>>
\dots @>{d_{q{-}1}}>> C^q(\mathfrak{g}, V) @>{d_q}>> \dots
\end{CD}
$$
where the differential $d_q$ is defined by:

\begin{equation}
\begin{split}
(d_q f)(X_1, \dots, X_{q{+}1})=
\sum_{i{=}1}^{q{+}1}(-1)^{i{+}1}
\rho(X_i)(f(X_1, \dots, \hat X_i, \dots, X_{q{+}1}))+\\
+ \sum_{1{\le}i{<}j{\le}q{+}1}(-1)^{i{+}j{-}1}
f([X_i,X_j],X_1, \dots, \hat X_i, \dots, \hat X_j, \dots, X_{q{+}1}).
\end{split}
\end{equation}

The cohomology of the complex $(C^*(\mathfrak{g}, V), d)$ is called
the cohomology of the Lie algebra $\mathfrak{g}$
with coefficients in the representation $\rho: \mathfrak{g} \to V$.

The cohomology of $(C^*(\mathfrak{g}, {\mathbb K}), d)$
($V= {\mathbb K}$ and $\rho: \mathfrak{g} \to {\mathbb K}$ is trivial)
is called the
cohomology with trivial coefficients of the Lie algebra
$\mathfrak{g}$ and is denoted by $H^*(\mathfrak{g})$.

One can remark that
$d_1: C^1(\mathfrak{g}, {\mathbb K}) \rightarrow
C^2(\mathfrak{g}, {\mathbb K})$ of the $(C^*(\mathfrak{g}, {\mathbb K}), d)$
is the dual mapping to the Lie bracket
$[ \, , ]: \Lambda^2 \mathfrak{g} \to \mathfrak{g}$.
Moreover the condition $d^2=0$ is equivalent to the Jacobi identity for $[,]$.

\begin{definition}
A Lie algebra $\mathfrak{g}$ is called $\mathbb N$-graded,
if it is decomposed to the direct sum of subspaces such that
$$\mathfrak{g}=\oplus_{i} \mathfrak{g}_{i}, \; i \in
{\mathbb N},
\quad \quad [\mathfrak{g}_{i}, \mathfrak{g}_{j}]
\subset \mathfrak{g}_{i+j}, \;
\forall \: i, j \in {\mathbb N}.
$$
\end{definition}
\begin{example}
The Lie algebra $\mathfrak{m}_0$ is defined
by its infinite basis $e_1, e_2, \dots, e_n, \dots $
with commutator relations:
$$ \label{m_0}
[e_1,e_i]=e_{i+1}, \; \forall \; i \ge 2.$$
The algebra $\mathfrak{m}_0$ has an infinite number of $\mathbb N$-gradings:
$$
\mathfrak{m}_0=\oplus_{i\in {\mathbb N}} {\mathfrak{m}_0}_i, \;
{\mathfrak{m}_0}_1=Span(e_1, e_2,\dots,e_k),\;
{\mathfrak{m}_0}_i=Span(e_{i{+}k{-}1}),
\:i \ge 2.
$$

\end{example}
\begin{remark}
We always omit the trivial commutator
relations $[e_i,e_j]=0$ in the definitions of Lie algebras.
\end{remark}

\begin{example}
Let us recall that the Witt algebra $W$ is
spanned by differential operators on the real line ${\mathbb R}^1$
with a fixed coordinate $x$
$$
e_i=x^{i+1}\frac{d}{dx}, \; i \in {\mathbb Z}, \quad
[e_i,e_j]= (j-i)e_{i{+}j}, \; \forall \;i,j \in {\mathbb Z}.
$$
We denote by $L_1$ a positive part of the Witt algebra, i.e.
$L_1$ is a subalgebra of $W$ spanned by all $e_i, \; i \ge 1$.

Obviously $W$ is a ${\mathbb Z}$-graded Lie algebra
with one-dimensional homogeneous
components:
$$
W=\oplus_{i\in {\mathbb Z}}W_i, \; W_i=Span(e_i).
$$
Thus $L_1$ is a ${\mathbb N}$-graded Lie algebra.
\end{example}

Let $\mathfrak{g}$ be a Lie algebra.
The ideals $C^{k}\mathfrak{g}$ of the descending central sequence define
a decreasing filtration $C$ of the Lie algebra
$\mathfrak{g}$
$$C^{1}\mathfrak{g}=\mathfrak{g} \supset C^{2}\mathfrak{g} \supset
\dots \supset C^{k}\mathfrak{g} \supset \dots; \qquad
[C^{k}\mathfrak{g},C^{l}\mathfrak{g}] \subset C^{k+l}\mathfrak{g}.$$
One can consider the associated graded Lie algebra
$${\rm gr}_C\mathfrak{g}=\bigoplus_{k=1} ({\rm gr}_C\mathfrak{g})_k, \;\;
({\rm gr}_C\mathfrak{g})_k=C^{k}\mathfrak{g}/C^{k{+}1}\mathfrak{g}.$$

\begin{proposition}
We have the following isomorphisms:
$$ {\rm gr}_C L_1 \cong {\rm gr}_C\mathfrak{m}_0 \cong
\mathfrak{m}_0.$$
\end{proposition}
\begin{remark}
$$
({\rm gr}_C\mathfrak{m}_0)_1=Span(e_1,e_2),\;
({\rm gr}_C\mathfrak{m}_0)_i=Span(e_{i{+}1}), \;i\ge 2.
$$
\end{remark}

Let $\mathfrak{g}=\oplus_{\alpha}\mathfrak{g}_{\alpha}$ be a
${\mathbb Z}$-graded Lie algebra
and $V=\oplus_{\beta} V_{\beta}$ is a ${\mathbb Z}$-graded
$\mathfrak{g}$-module, i.e.,
$\mathfrak{g}_{\alpha}V_{\beta} \subset V_{\alpha {+} \beta}.$
Then the complex $(C^*(\mathfrak{g}, V), d)$
can be equipped with the ${\mathbb Z}$-grading
$C^q(\mathfrak{g},V) =
\bigoplus_{\mu} C^q_{(\mu)}(\mathfrak{g},V)$, where
a $V$-valued $q$-form $c$ belongs to
$C^q_{(\mu)}(\mathfrak{g},V)$
iff for
$X_1 \in \mathfrak{g}_{\alpha_1},
\dots,  X_q \in \mathfrak{g}_{\alpha_q}$ we have
$$c(X_1,\dots,X_q) \in
V_{\alpha_1{+}\alpha_2{+}\dots{+}\alpha_q{+}\mu}.$$

This grading is compatible with the differential  $d$ and
hence we have ${\mathbb Z}$-grading in the cohomology:
$$
H^{q} (\mathfrak{g},V)= \bigoplus_{\mu \in {\mathbb Z}}
H^{q}_{(\mu)} (\mathfrak{g},V).
$$

\begin{remark}
The trivial $\mathfrak{g}$-module ${\mathbb K}$ has only one non-trivial
homogeneous component ${\mathbb K}={\mathbb K}_0$.
\end{remark}

The exterior product in $\Lambda^*(\mathfrak{g})$
induces a structure of a
bigraded algebra in the cohomology
$H^*(\mathfrak{g})$:
$$
H^{q}_{k} (\mathfrak{g}) \wedge
H^{p}_{l} (\mathfrak{g}) \to
H^{q{+}p}_{k+l} (\mathfrak{g}).
$$

Let $\mathfrak{g}=
\oplus_{\alpha >0} \mathfrak{g}_{\alpha}$
be a ${\mathbb N}$-graded Lie algebra
and $V=\oplus_{\beta} V_{\beta}$ is a ${\mathbb Z}$-graded
$\mathfrak{g}$-module.

One can define
a decreasing filtration ${\mathcal F}$ of
$(C^*(\mathfrak{g},V),d)$:
$$
{\mathcal F}^0 C^*(\mathfrak{g},V)\supset \dots
\supset {\mathcal F}^q C^*(\mathfrak{g},V)
\supset {\mathcal F}^{q{+}1} C^*(\mathfrak{g},V) \supset \dots
$$
where the subspace ${\mathcal F}^q C^{p{+}q}(\mathfrak{g},V)$
is spanned by $p{+}q$-forms $c$ in $C^{p{+}q}(\mathfrak{g},V)$
such that
$$
c(X_1,\dots,X_{p{+}q}) \in \bigoplus_{\alpha=q} V_{\alpha}, \;
\forall X_1,\dots,X_{p{+}q} \in \mathfrak{g}.
$$

The filtration ${\mathcal F}$ is compatible with $d$.

Let us consider the corresponding spectral sequence $E_r^{p,q}$:
\begin{proposition}[\cite{Fu},\cite{FeFuRe}]
$E_1^{p,q}=V_q \otimes H^{p{+}q}({\mathfrak g})$.
\end{proposition}
\begin{proof}
We have the following natural isomorphisms:
\begin{equation}
\begin{split}
C^{p{+}q}({\mathfrak g},V)
= V \otimes \Lambda^{p{+}q}({\mathfrak g}^*)\\
E_0^{p,q}={\mathcal F}^q C^{p{+}q}(\mathfrak{g},V)/
{\mathcal F}^{q{+}1} C^{p{+}q}(\mathfrak{g},V)
=V_q \otimes \Lambda^{p{+}q}({\mathfrak g}^*).
\end{split}
\end{equation}

Now the proof follows from the formula for the
$d_0^{p,q}: E_0^{p,q} \to E_0^{p{+}1,q}$:
$$d_0(v \otimes f)=v\otimes df,$$
where $v \in V, f \in \Lambda^{p{+}q}({\mathfrak g}^*)$
and $df$ is the standard
differential of the cochain complex of $\mathfrak{g}$ with trivial
coefficients.
\end{proof}

\section{Cohomology $H^*(L_1)$ and $H^*(\mathfrak{m}_0)$}
\label{L_1_m_0}
\begin{theorem}[Goncharova,\cite{G}]
The Betti numbers ${\rm dim} H^q(L_1)=2$, for every
$q \ge 1$, more precisely
$$ {\rm dim} H_k^q(L_1)=
\left\{\begin{array}{r}
   1, \hspace{0.6em}{\rm if}~k=\frac{3q^2 \pm q}{2}, \\
   0,\hspace{1.76em}{\rm otherwise.}\\
   \end{array} \right . \hspace{3.3em} $$
\end{theorem}
We will denote in the sequel by $g^q_{\pm}$ the generators of the spaces
$H^q_{\frac{3q^2\pm q}{2}}(L_1)$.
The numbers
$\frac{3q^2 \pm q}{2}$ are so called Euler pentagonal numbers.
A sum of two arbitrary pentagonal numbers is not a pentagonal number,
hence the algebra
$H^*(L_1)$ has a trivial multiplication.

\begin{example}

1) $H^1(L_1)$ is generated
by $g^1_-=[e^1]$ and $g^1_+=[e^2]$;

2) the basis of $H^2(L_1)$ consists of two
classes $g^2_-=[e^1 {\wedge} e^4]$ and
$g^2_+=[e^2 {\wedge} e^5-3e^3 {\wedge} e^4]$ of weights $5$ and $7$
respectively.
\end{example}

The cohomology algebra
$H^*(\mathfrak{m}_0)$ was calculated by Fialowski and Millionschikov
in \cite{FialMill}.

It were introduced two operators in \cite{FialMill}:

1) $D_1: \Lambda^*(e_2, e_3, \dots)
\to \Lambda^*(e_2, e_3, \dots)$,
\begin{equation}
\begin{split}
D_1(e^2)=0, \; D_1(e^i)= e^{i-1}, \; \forall i\ge 3,\\
D_1(\xi\wedge \eta)=D_1(\xi)\wedge \eta +\xi\wedge D_1(\eta),
\; \: \forall \xi, \eta \in \Lambda^*(e_2, e_3, \dots).
\end{split}
\end{equation}

2) and its right inverse $D_{-1}: \Lambda^*(e^2,e^3,\dots) \to
\Lambda^*(e^2,e^3,\dots)$,
\begin{equation}
\begin{split}
\label{D_{-1}}
D_{-1}e^i=e^{i+1}, \;
D_{-1}(\xi{\wedge} e^i)=
\sum_{l\ge 0}(-1)^l D_{1}^{l}(\xi){\wedge} e^{i+1+l},
\end{split}
\end{equation}
where $i\ge 2$ and $\xi$ is an arbitrary form in
$\Lambda^*(e^2,\dots,e^{i-1})$.
The sum in the definition (\ref{D_{-1}}) of $D_{-1}$ is always finite
because $D_1^l$ decreases the second grading by $l$.
For instance,
$$D_{-1}(e^i\wedge e^k)=\sum_{l=0}^{i-2} ({-}1)^l e^{i-l}{\wedge} e^{k+l+1}.$$

\begin{proposition}
The operators $D_1$ and $D_{-1}$ have the following properties:
$$
\label{D_1,D_{-1}}
d\xi=e^1\wedge D_1\xi, \; e^1\wedge \xi=d D_{-1}\xi, \; D_1D_{-1}\xi=\xi,
\quad \xi \in \Lambda^*(e^2,e^3,\dots).
$$
\end{proposition}

\begin{theorem}[\cite{FialMill}]
The infinite dimensional bigraded cohomology
$H^*(\mathfrak{m}_0)=\oplus_{k,q} H^q_k(\mathfrak{m}_0)$
is spanned by the cohomology classes of
$e^1$, $e^2$ and of the following homogeneous cocycles:
\begin{equation}
\label{cocycles}
\omega(e^{i_1}{\wedge}\dots {\wedge} e^{i_q} {\wedge}
e^{i_q{+}1})=
\sum\limits_{l\ge 0}(-1)^l D_1^l(e^{i_1}\wedge \dots
\wedge e^{i_q})\wedge e^{i_q+1+l},
\end{equation}
where $q \ge 1, \; 2\le i_1 {<}i_2{<}{\dots} {<}i_q$.
\end{theorem}

Formula (\ref{cocycles})
determines a homogeneous closed $(q{+}1)$-form
of the second grading $i_1{+}{\dots}{+}i_{q{-}1}{+}2i_{q}{+}1$.
It has only one monomial in its decomposition
of the form $\xi \wedge e^i \wedge e^{i+1}$ and it is
$e^{i_1}{\wedge}\dots {\wedge} e^{i_q} {\wedge}
e^{i_q{+}1}$.

The whole number of linearly indepenedent
$q$-cocycles of the second grading $k{+}\frac{q(q{+}1)}{2}$ is equal to
$${\rm dim} H_{k{+}\frac{q(q{+}1)}{2}}^q(\mathfrak{m}_0)=
P_q(k)-P_q(k{-}1),$$
where $P_q(k)$ denotes the number of
(unordered) partitions of a positive integer $k$ into $q$ parts.

\begin{example}
\label{3cocycle}
\begin{multline*}
\omega(e^5{\wedge} e^6 {\wedge} e^7)=
e^5{\wedge} e^6 {\wedge}e^7-e^4{\wedge} e^6 {\wedge} e^8+(e^3{\wedge} e^6+
e^4 {\wedge} e^5){\wedge} e^9-\\
-(e^2{\wedge}e^6+2e^3 {\wedge}e^5){\wedge}e^{10}
+(3e^2{\wedge} e^5+2e^3{\wedge} e^4){\wedge} e^{11}-
5e^2{\wedge} e^4{\wedge} e^{12}+5e^2{\wedge} e^3{\wedge} e^{13}.
\end{multline*}
\end{example}

The multiplicative structure in $H^*(\mathfrak{m}_0)$
was also found in \cite{FialMill} explicitly. In particular
$$
[e^1] {\wedge} \omega(\xi{\wedge}e^i{\wedge}e^{i+1})=0,\;
[e^2] {\wedge} \omega(\xi{\wedge}e^i{\wedge}e^{i+1})
=\omega(e^2{\wedge}\xi{\wedge}e^i{\wedge}e^{i+1}).
$$

\section{Massey products in cohomology.}
In this section we follow \cite{JPM} and \cite{BaTa} presenting the
definitions of Massey products.
\label{Massey_S}
Let $\A=\oplus_{l \ge 0}\A^l$
be a differential graded algebra over a field ${\mathbb K}$.
It means that the following operations are defined:
an associative multiplication
$$
\wedge:\A^l\times\A^m\to\A^{l+m},\; l,m\geq 0, \; l,n \in {\mathbb Z}.
$$
such that $a\wedge b=(-1)^{lm}\,b\wedge a$ for $a\in\A^l$, $b\in\A^m$,
and a differential $d, \; d^2=0$
$$
d:\A^l\to\A^{l+1},\ \ l\geq 0,
$$
satisfying the Leibniz rule
$d\,(a\wedge b)=d\,a\wedge b+(-1)^l a\wedge d\,b$ for $a\in\A^l$.
\begin{example}
$\A=\Lambda^*(\mathfrak{g})$ is the cochain complex of a Lie algebra.
\end{example}

For a given differential graded algebra  $(\A,d)$ we denote
by $T_n(\A)$ a space of all upper triangular
$(n+1)\times (n+1)$-matrices with entries from $\A$,
vanishing at the main diagonal. $T_n(\A)$ has a structure of
a differential algebra with a standard matrix multiplication,
where matrix entries are multiplying as elements of $\A$.
A differential $d$ on $T_n(\A)$ is defined by
\begin{equation}
d\,A=(d\,a_{ij})_{1\leq i,j\leq n{+}1}.
\label{1.2.1}
\end{equation}

An involution $a\to\bar{a}=(-1)^{k+1}a, a \in A^k$
of $\A$ can be extended to an involution of $T_n(\A)$ as
$\bar{A}=(\bar{a}_{ij})_{1\leq i,j\leq n{+}1}$.
It satisfies the following properties:
$$
\overline{\bar{A}}=A, \quad
\overline{AB}=-\bar{A}\bar{B}, \quad \overline{d\,A}=-d\,\bar{A}.
$$
Also we have the generalized Leibniz rule for the differential
(\ref{1.2.1})
$$
d\,(AB)=(d\,A)B-\bar{A}(d\,B).
$$

The algebra $T_n(\A)$ has a
two-sided center $I_n(\A)$ of matrices
$$
\left(\begin{array}{cccc}
0 & \dots & 0 & \tau \\
0 & \dots & 0 & 0 \\
& \dots& & \\
0 & \dots & 0 & 0
\end{array}\right), \quad \tau \in \A
$$

\begin{definition}[\cite{BaTa}]
A matrix $A \in T_n(\A)$ is called the matrix of
a formal connection if it satisfies the Maurer-Cartan equation
\begin{equation}
\mu(A)=d\,A-\bar{A}\cdot A \in I_n(\A).
\label{star}
\end{equation}
\end{definition}
\begin{proposition}[\cite{BaTa}]
Let $A$ be the matrix of a formal connection,
then the entry $\tau \in \A$ of the matrix $\mu(A) \in I_n(\A)$
in the definition (\ref{star}) is closed.
\end{proposition}
\begin{proof}
We have the following generalized Bianci identity
for the Maurer-Cartan operator
$\mu(A)=d\,A - \bar{A} \cdot A$ ($A$ is an arbitrary matrix):
$$
d\,\mu(A)=\overline{\mu(A)}\cdot A+A\cdot\mu(A).
$$
Indeed it's easy to verify the following equalities:
$$
d\,\mu(A)=-d\,(\bar{A}\cdot A)=
-d\,\bar{A}\cdot A+A\cdot d\,A=\overline{d\,A}\cdot A+A\cdot d\,A=
$$
$$
=\overline{(\mu(A)+\bar{A}\cdot A)}\cdot A+A(\mu(A)+\bar{A}\cdot A)=
$$
$$
= \overline{\mu(A)}\cdot A-
A\cdot\bar{A}\cdot A+A\cdot\mu{A}+A\cdot\bar{A}\cdot A=
$$
$$
=\overline{\mu(A)}\cdot A+A\cdot\mu(A).
$$
\end{proof}
Now let $A$ be the matrix of a formal connection, then
the matrix $\mu(A)$ belongs to the center $I_n(\A)$ and hence $d\mu(A)=0$.
One can think of $\mu(A)$ as the curvature
matrix of a formal connection $A$.

Let $A$ be an upper triangular matrix from $T_n(\A)$.
One can rewrite it in the following notation:
$$
A=\left(\begin{array}{cccccc}
  0      & a(1,1) & a(1,2) & \dots & a(1,n-1)   & a(1,n)   \\
  0      & 0      & a(2,2) & \dots & a(2,n-1)   & a(2,n)   \\
  \dots  & \dots  & \dots  & \dots & \dots      & \dots    \\
  0      & 0      & 0      & \dots & a(n-1,n-1) & a(n-1,n) \\
  0      & 0      & 0      & \dots & 0          & a(n,n)   \\
  0      & 0      & 0      & \dots & 0          & 0
  \end{array}\right).
$$

\begin{proposition}
A matrix $A \in T_n(\A)$ is the matrix of a formal connection if and only if
the following conditions on its entries hold on
\begin{equation}
\label{def_syst}
\begin{split}
a(i,i)=a_i \in \A^{p_i}, \quad i=1,\dots,n;\\
a(i,j)\in\A^{p(i,j)+1}, \quad p(i,j)=\sum^j_{r=i}(p_r-1);\\
d\,a(i,j)=\sum_{r=i}^{j-1}\bar{a}(i,r)\cdot a(r+1,j),\;\;
(i,j)\neq(1,n).
\end{split}
\end{equation}
\end{proposition}
\begin{proof}
The system (\ref{def_syst}) is  just the Maurer-Cartan equation
rewritten in terms of the entries of the matrix $A$ and it
is a part of the classical definition \cite{K} of
the defining system for a Massey product.
\end{proof}
\begin{definition}[\cite{K}]
A collection of elements,
$A=(a(i,j))$, for $1\leq i\leq j\leq n$ and $(i,j)\neq(1,n)$
is said to be a defining  system for the product
$\langle a_1,\dots,a_n\rangle$ if it satisfies (\ref{def_syst}).

In this situation
the $(p(1,n)+2)$-dimensional cocycle
$$
c(A)=\sum_{r=1}^{n-1}\bar{a}(1,r)a(r+1,n)
$$
is called the related cocycle of the defining system $A$.
\end{definition}
\begin{remark}
We saw that the notion of the defining system is equivalent to
the notion of the formal connection. However one has to remark
that an entry $a(1,n)$ of the matrix $A$ of a formal connection
does not belong to the corresponding defining system $A$,
it can be taken as an arbitrary element from $\A$.
In this event the only one (possible) nonzero entry
$\tau$ of the Maurer-Cartan matrix $\mu(A)$ is equal to $-c(A)+da(1,n)$.
\end{remark}

\begin{definition}[\cite{K}]
The $n$-fold product $\langle a_1,\dots,a_n\rangle$ is defined if
there is at least one defining system for it
(a formal connection $A$ with
entries $a_1,\dots,a_n$ at the second diagonal). If it is defined, then
$\langle a_1,\dots,a_n\rangle$ consists of all cohomology
classes $\alpha \in H^{p(1,n)+2}(\A)$ for which
there exists a defining system $A$
(a formal connection $A$)
such that $c(A)$ ($-\tau$ respectively) represents $\alpha$.
\end{definition}

\begin{theorem}[\cite{K},\cite{BaTa}]
The operation $\langle a_1,\dots,a_n\rangle$ depends only on the
cohomology classes of the elements $a_1,\dots,a_n$.
\end{theorem}
\begin{proof}
A changing of an arbitrary entry $a_{ij}, \:j > i$
of the matrix $A$ of a formal connection to $a_{ij}+db$ leads
to a replacement of $A$ by
$$
A'=A+db\cdot E_{ij}+A\cdot b\cdot E_{ij}-\bar b\cdot E_{ij}\cdot A,
$$
where $E_{ij}$ is a scalar matrix wich has $1$
on $(i,j)$-th place and zeroes on all others.
For the corresponding Maurer-Cartan matrix we will have
$$
\mu(A')=\mu(A)+d((A\cdot b \cdot E_{ij}-\bar b \cdot E_{ij}\cdot A)\cap I_n).
$$
\end{proof}
\begin{definition}[\cite{K}]
A set of closed elements $a_i, i=1,\dots,n$ from $\A$ representing
some cohomology classes ${\alpha}_i \in H^{p_i}(\A), i=1,\dots,n$
is said to be a defining system for
the Massey $n$-fold product $\langle {\alpha}_1,\dots,{\alpha}_n\rangle$
if it is one for $\langle a_1,\dots,a_n\rangle$. The Massey
$n$-fold product $\langle {\alpha}_1,\dots,{\alpha}_n\rangle$
is defined if $\langle a_1,\dots,a_n\rangle$ is defined, in which case
$\langle {\alpha}_1,\dots,{\alpha}_n\rangle=\langle a_1,\dots,a_n\rangle$
as subsets in $H^{p(1,n)+2}(\A)$.
\end{definition}

\begin{example}
For $n=2$ the matrix $A$ of a formal connection has a form
$A=\left(\begin{array}{ccc}
0 & a & c\\
0 & 0 & b\\
0& 0& 0
\end{array}\right)$ and the matrix Maurer-Cartan equation is equivalent
to two equations $da=0$ and $db=0$. Evidently
$\langle \alpha, \beta \rangle = \bar \alpha \cdot \beta$.
\end{example}

\begin{example}[Triple Massey products]
Let $\alpha$, $\beta$, and  $\gamma$ be the cohomology classes of closed
elements $a \in\A^p$, $b\in\A^q$, and  $c\in\A^r$.
The Maurer-Cartan equation for
$$
A=\left(\begin{array}{cccc}
0 & a & f & h \\
0 & 0 & b & g \\
0 & 0 & 0 & c \\
0 & 0 & 0 & 0
\end{array}\right).
$$
is equivalent to
\begin{equation}
d\,f=(-1)^{p+1}\,a\wedge b,\ \ d\,g=(-1)^{q+1}\,b\wedge c.
\label{ast}
\end{equation}
Hence the triple Massey product
$\langle\alpha,\beta,\gamma\rangle$
is defined if and only if
$$
\alpha\cdot \beta=\beta \cdot \gamma=0\ \ \mbox{in}\ \ H^{\ast}(\A).
$$
If these conditions are satisfied then
the Massey product $\langle \alpha, \beta, \gamma \rangle$ is defined
as a subset in $H^{p{+}q{+}{r}{-}1}(\A)$ of the following form
$$
\langle \alpha, \beta, \gamma \rangle=\left\{
[(-1)^{p+1} a\wedge g+(-1)^{p+q} f\wedge c]\right\}.
$$
Since $f$ and  $g$ are defined by (\ref{ast}) up to closed elements
from $\A^{p+q-1}$  and $\A^{q+r-1}$ respectively, the triple Massey product
$\langle \alpha,\beta,\gamma\rangle$ is an affine subspace of
$H^{p{+}q{+}{r}{-}1}(\A)$ parallel to
$\alpha \cdot H^{q+r-1}(\A)+ H^{p+q-1}(\A)\cdot \gamma$.
\end{example}
\begin{remark}
We defined Massey products as the multi-valued operations in general.
Moreoften in the literature the triple Massey product is defined
as a quotient-space
$\langle \alpha, \beta, \gamma \rangle/
(\alpha \cdot H^{q+r-1}(\A)+ H^{p+q-1}(\A)\cdot \gamma)$ and it is
single-valued in this case (see \cite{Fu}).
\end{remark}

\begin{definition}
Let an  $n$-fold Massey product
$\langle {\alpha}_1,\dots,{\alpha}_n\rangle$
be defined. It is called trivial if it contains the trivial cohomology class:
$0\in\langle {\alpha}_1,\dots,{\alpha}_n\rangle$.
\end{definition}

\begin{proposition}
\label{triviality}
Let a Massey product
$\langle {\alpha}_1,\dots,{\alpha}_n\rangle$
is defined. Then all Massey products
$\langle {\alpha}_l,\dots,{\alpha}_q\rangle, 1\le l < q \le n, q-l<n-1$
are defined and trivial.
\end{proposition}
\begin{proof}
It follows from (\ref{def_syst}).
\end{proof}
\begin{remark}
The triviality of all Massey products
$\langle {\alpha}_l,\dots,{\alpha}_q\rangle, 1\le l < q \le n, q-l<n-1$
is only a necessary condition
for a Massey product $\langle {\alpha}_1,\dots,{\alpha}_n\rangle$
to be defined. It is sufficient only in the case $n=3$.
\end{remark}

Let us denote by $GT_n({\mathbb K})$ a group of non-degenerate
upper triangular $(n{+}1,n{+}1)$-matrices of the form:
$$
C=\left(\begin{array}{ccccc}
  a_{1,1}      & a_{1,2} & \dots & a_{1,n}   & a_{1,n{+}1}   \\
  0      & a_{2,2}      &  \dots & a_{2,n}  & a_{2,n{+}1}   \\
 \dots   &        &       \dots &            & \dots     \\
  0      & 0            & \dots & a_{n,n}& a_{n,n{+}1}   \\
  0      & 0            & \dots & 0          & a_{n{+}1,n{+}1}
  \end{array}\right).
$$
\begin{proposition}
Let $A \in T_n({\A})$ be the matrix of a formal connection and
$C$ an arbitrary matrix from $GT_n({\mathbb K})$. Then the matrix
$C^{-1}AC\in T_n({\A})$ and satisfies the Maurer-Cartan equation,
i.e. is again the matrix of a formal connection.
\end{proposition}
\begin{proof}
$$
d(C^{-1}AC)-\bar C^{-1}\bar A\bar C \wedge C^{-1}AC=
C^{-1}\left(dA-\bar A\wedge A\right)C=0.
$$
\end{proof}
\begin{example}
Let $A \in T_n({\A})$ be the matrix of a formal connection
(defining system) for a Massey product
$\langle \alpha_1,\dots, \alpha_n\rangle$.
Then a matrix $C^{-1}AC$ with
$$
C=\left(\begin{array}{ccccc}
  1      & 0 & \dots & 0   & 0   \\
  0      & x_1      &  \dots & 0  & 0   \\
 \dots   &        &       \dots &            & \dots     \\
  0      & 0            & \dots &  x_1{\dots}x_{n{-}1}& 0   \\
  0      & 0            & \dots & 0          & x_1{\dots}x_{n{-}1}x_n
  \end{array}\right)
$$
is a defining system for $\langle x_1\alpha_1,\dots, x_n\alpha_n\rangle=
x_1\dots x_n\langle x_1\alpha_1,\dots, x_n\alpha_n\rangle$.
\end{example}

\begin{definition}
\label{A_equiv}
Two matrices $A$ and $A'$ of formal connections
from  $T_n({\A})$ are equivalent if there exists a
matrix $C \in GL(n{+}1,{\mathbb K})$ such that
$$
A'=C^{-1}AC.
$$
\end{definition}

\begin{example}
Triple products $\langle\alpha,\beta,\gamma\rangle$ and
$\langle x \alpha,y \beta,z \gamma\rangle$, where
$x,y,z \ne 0$, are equivalent with
$$
C=\left(\begin{array}{cccc}
  1      & 0 &  0   & 0   \\
  0      & x      &   0  & 0   \\
  0      & 0            &  xy & 0   \\
  0      & 0            &  0          & xyz
  \end{array}\right)
$$
and
$$
\langle x \alpha,y \beta,z \gamma\rangle=
xyz\langle \alpha, \beta, \gamma \rangle, \quad x,y,z \in{\mathbb K}.
$$
\end{example}
\begin{remark}
Following the original Massey work \cite{Mass}
some higher order cohomological operations
that we call now Massey products were introduced
in the 60s in \cite{K} and \cite{JPM}.
The relation between Massey products and the
Maurer-Cartan equation was first noticed by May \cite{JPM} and this analogy
was not developed untill \cite{BaTa}.

In the present article
we deal only with Massey products of non-trivial cohomology classes.
It is possible to take some of them trivial, but in this
situation is more natural to work with so-called matrix
Massey products that were first introduced by May \cite{JPM}. This
approach was also developped in \cite{BaTa}. We will not treat this
case in the sequel.
\end{remark}

\section{Massey products and thread modules}
\label{g-modules}
Let $T_n({\mathbb K})$ be a
Lie algebra of upper triangular $(n+1,n+1)$-matrices over
a field $\mathbb K$ of zero characteristic and
$\rho: \mathfrak{g} \to T_n({\mathbb K})$ be a representation of
a Lie algebra $\mathfrak{g}$.

\begin{example}
We take $n=1$ and consider a linear map
$$
\rho: x \in \mathfrak{g} \to
\left(\begin{array}{cc}
0 & a(x) \\
0& 0
\end{array}\right).$$
It is evident that
$\rho$ is a Lie algebra homomorphism if and only if
the linear form $a \in {\mathfrak{g}}^*$ is closed
$$
da(x,y)=a([x,y])=a(x)a(y)-a(y)a(x)=0, \; \forall x,y \in \mathfrak{g}.
$$

In other words the matrix
$A=\left(\begin{array}{cc}
0 & a \\
0& 0
\end{array}\right)$
satisfies the "strong" Maurer-Cartan equation
$dA-\bar A \wedge A =0$.
\end{example}
\begin{remark}
We recall that we defined in the Section \ref{Massey_S}
the involution of a graded $\A$ as $\bar a=(-1)^{k{+}1} a, a \in \A^k$.
Thus for a matrix $A$ with entries from ${\mathfrak{g}}^*$
we have $\bar A =A$. One has to remark that $\bar a$ differs
by the sign from the definition of $\bar a$ in \cite{K},
however in \cite{JPM2} one meets our sign rule.
\end{remark}

\begin{proposition}
A matrix $A$ with entries from ${\mathfrak{g}}^*$ defines
a representation $\rho: \mathfrak{g} \to T_n({\mathbb K})$
if and only if $A$ satisfies the strong Maurer-Cartan equation
$$dA-\bar A \wedge A=0.$$
\end{proposition}
\begin{proof}
$$
(dA-\bar A \wedge A)(x,y)=A([x,y])-\left[ A(x), A(y)\right],\;\forall x,y \in
\mathfrak{g}.
$$
\end{proof}
\begin{example}
For $n=2$ the matrix $A$ of a representation
$\rho$ has a form
$A=\left(\begin{array}{ccc}
0 & a & c\\
0 & 0 & b\\
0& 0& 0
\end{array}\right),$ where $a,b,c \in {\mathfrak{g}}^*$
and the strong Maurer-Cartan equation is equivalent
to the following equations on entries $a,b,c$:
$$da=0, \quad db=0, \quad dc= a \wedge b.$$
\end{example}

The Lie algebra $T_n({\mathbb K})$ has a
one-dimensional center $I_n({\mathbb K})$ spanned by the matrix
$$
\left(\begin{array}{cccc}
0 & \dots & 0 & 1 \\
0 & \dots & 0 & 0 \\
& \dots& & \\
0 & \dots & 0 & 0
\end{array}\right).
$$
One can consider an one-dimensional central extension
$$
\begin{CD}
0 @>>> \mathbb K \cong I_n({\mathbb K}) @>>>
T_n(\mathbb K) @>{\pi}>> \tilde T_n(\mathbb K) @>>> 0.
\end{CD}
$$

\begin{proposition}[\cite{FeFuRe}, \cite{Dw}]
Fixing a Lie algebra homomorphism $\tilde \varphi : \mathfrak{g} \to
\tilde T_n({\mathbb K})$ is equivalent to fixing a defining system
$A$ with elements from $\mathfrak{g}^*=\Lambda^1(\mathfrak{g})$.
The related cocycle $c(A)$ is cohomologious to zero
if and only if $\tilde \varphi$ can be lifted to a
homomorphism $\varphi : \mathfrak{g} \to  T_n({\mathbb K})$,
$\tilde \varphi = \pi \varphi$.
\end{proposition}

Taking a $(n{+}1)$-dimensional linear space $V$ over ${\mathbb K}$
and a representation $\varphi : \mathfrak{g} \to T_n({\mathbb K})$
one gets a $\mathfrak{g}$-module structure of $V$
defined in the coordinates $x=(x_1,\dots,x_{n{+}1})$ with respect to
some fixed basis $v_1,\dots,v_{n{+}1}$ of $V$
$$
gv=\varphi(g)x,\quad g \in \mathfrak{g}, \; v \in V.
$$

We recall the standard definition from the representations theory.
\begin{definition}
Two representations
$\varphi : \mathfrak{g} \to T_n({\mathbb K})$ and
$\varphi' : \mathfrak{g} \to T_n({\mathbb K})$ are called isomorphic
(equivalent) if there exist two bases $v_1,\dots,v_{n{+}1}$ and
$v_1',\dots,v_{n{+}1}'$ in a $(n{+}1)$-dimensional linear space $V$
such that the corresponding $\mathfrak{g}$-module structures
coinside $\rho(g)v=\rho'(g)v$. Or in other words, if there exists
a matrix $C \in GL(n{+}1,{\mathbb K})$ such that
$$
\varphi'(g)=C^{-1}\varphi(g)C, \; \forall g \in \mathfrak{g}.
$$
\end{definition}
It is evident that this definition
is equivalent to the Definition \ref{A_equiv} when we consider
a Massey product $\langle \omega_1,\dots, \omega_n \rangle$
of $1$-cohomology classes $\omega_1,\dots, \omega_n$.

\begin{proposition}
$$
\left[{\rm gr}_C\mathfrak{g},{\rm gr}_C\mathfrak{g}\right]=
\oplus_{i\ge 2}({\rm gr}_C\mathfrak{g})_i.
$$
\end{proposition}
\begin{proof}
An inclusion
$
\left[{\rm gr}_C\mathfrak{g},{\rm gr}_C\mathfrak{g}\right] \subset
\oplus_{i\ge 2}({\rm gr}_C\mathfrak{g})_i
$
is evident for the graded Lie algebra ${\rm gr}_C\mathfrak{g}=
\oplus_{i\ge 1}({\rm gr}_C\mathfrak{g})_i$.
Hence it is sufficient to prove an inclusion
$
({\rm gr}_C\mathfrak{g})_i \subset
\left[{\rm gr}_C\mathfrak{g},{\rm gr}_C\mathfrak{g}\right]
$
for an arbitrary $i \ge 2$.
The last one in its turn follows from
$$
({\rm gr}_C\mathfrak{g})_i=
[\mathfrak{g},C^{i{-}1}\mathfrak{g}]+C^{i{+}1}\mathfrak{g}=
[\mathfrak{g}{+}C^{2}\mathfrak{g},C^{i{-}1}\mathfrak{g}{+}C^{i}\mathfrak{g}].
$$
\end{proof}

\begin{corollary}
We have isomorphisms
$$
\label{h_1_gr}
H^1(\mathfrak{g}) \cong
\left(\mathfrak{g}/[\mathfrak{g},\mathfrak{g}]\right)^* \cong
({\rm gr}_C\mathfrak{g})_1^* \cong H^1({\rm gr}_C\mathfrak{g}).
$$
\end{corollary}

\begin{definition}
A homomorphism  $\rho: \mathfrak{g} \to \mathfrak{h}$
of two ${\mathbb N}$-graded Lie algebras
$\mathfrak{g}=\oplus_i \mathfrak{g}_i$ and
$\mathfrak{h}=\oplus_i \mathfrak{h}_i$
is called graded if
$$
\rho (\mathfrak{g}_i) \subset \mathfrak{h}_i, \quad \forall i \in {\mathbb N}.
$$
\end{definition}

\begin{definition}
Let $\rho: \mathfrak{g} \to T_n({\mathbb K})$ be a representation.
We call a representation of graded Lie algebras
$
\tilde \rho :{\rm gr}_C\mathfrak{g} \to {\rm gr}_C T_n({\mathbb K}),
$
defined by the rule
$$
\tilde \rho (x{+}C^{k{+}1}\mathfrak{g})=\rho(x){+}C^{k{+}1}T_n({\mathbb K}),
\; x {+} C^{k{+}1}\mathfrak{g} \in
C^{k}\mathfrak{g}/C^{k{+}1}\mathfrak{g}, \; k \ge 1,
$$
an associated graded representation to $\rho$.
\end{definition}

\begin{remark}
Let $\tilde \rho: \mathfrak{g} \to T_n({\mathbb K})$ be some
representation of a Lie algebra $\mathfrak{g}$
such as $\mathfrak{g}\cong {\rm gr}_C\mathfrak{g}$. It is not hard
to describe the corresponding matrix $A$ of a formal connection:

1) the first diagonal is of zeroes (like all matrices from $T_n({\mathbb K})$);

2) the second one contains only elements from $\mathfrak{g}^*$
of degree one;

3) the $k$-th diagonal consists only  of elements of degree $k-1$.
\end{remark}

\begin{example}
We take $n=3$ and a homomorphism $\rho: \mathfrak{m}_0 \to T_3({\mathbb K})$
is defined by
$$
\rho(e_1)=
\left(\begin{array}{cccc}
0 & 0 & 0& 1\\
0 & 0 & 1& 0\\
0 & 0 & 0& 1\\
0 & 0 & 0& 0 \\
\end{array}\right), \;
\rho(e_2)=
\left(\begin{array}{cccc}
0 & 1 & 0& 0\\
0 & 0 & 0& 1\\
0 & 0 & 0& 0\\
0 & 0 & 0& 0 \\
\end{array}\right).
$$
On can to define $\rho$ only on $e_1$ and $e_2$ because the algebra
$\mathfrak{m}_0$ is generated by these two elements.

The corresponding matrix $A$ of a formal connection is equal to
$$
A=\left(\begin{array}{cccc}
0 & e^2 & -e^3& e^4+e^1\\
0 & 0 & e^1& e^2\\
0 & 0 & 0& e^1\\
0 & 0 & 0& 0 \\
\end{array}\right).
$$
For the associated graded representation
$\tilde \rho: \mathfrak{m}_0 \to T_3({\mathbb K})$ we have
$$
\tilde A=\left(\begin{array}{cccc}
0 & e^2 & {-}e^3& e^4\\
0 & 0 & e^1& 0\\
0 & 0 & 0& e^1\\
0 & 0 & 0& 0 \\
\end{array}\right).
$$
\end{example}
We recall that elements $e^1$ and $e^2$ have grading one, $e^3$ and $e^4$
have gradings $2$ and $3$ in
$\Lambda^*(\mathfrak{m}_0)=\Lambda^*({\rm gr}_C\mathfrak{m}_0)$ respectively.

Identifying the spaces $H^1(\mathfrak{g})$ and
$H^1({\rm gr}_C\mathfrak{g})$ we come to the following
proposition:
\begin{proposition}
Let a Massey product $\langle \omega_1,\omega_2,{\dots},\omega_n \rangle$
be defined and trivial in $H^2(\mathfrak{g})$ for some $1$-cohomology classes
$\omega_i \in H^1(\mathfrak{g})$ of a Lie algebra $\mathfrak{g}$. Then
$\langle \omega_1,\omega_2,{\dots},\omega_n \rangle$
is also defined and trivial
in $H^2({\rm gr}_C\mathfrak{g})$.
\end{proposition}

\begin{proposition}
Let $\mathfrak{g}$ be a Lie algebra such that
$\mathfrak{g}\cong {\rm gr}_C\mathfrak{g}$
and a Massey product $\langle \omega_1,\omega_2,{\dots},\omega_n \rangle$
be defined and trivial for some $\omega_i \in H^1(\mathfrak{g})$.
Then there exists a graded defining system $A$ (the matrix $A$
of a graded formal connection) for
$\langle \omega_1,\omega_2,{\dots},\omega_n \rangle$.
\end{proposition}

\begin{definition}
A thread module over a ${\mathbb N}$-graded Lie algebra
$\mathfrak{g}=\oplus_i \mathfrak{g}_i$ is a
${\mathbb N}$-graded $\mathfrak{g}$-module $V=\oplus_{i\in {\mathbb N}} V_i$
such as
$$
\dim V_i=1,\quad
\mathfrak{g}_i V_j \subset V_{i+j},
\; \forall i,j \in {\mathbb N}.
$$
\end{definition}

Fixing a basis $\{f_j\}, \;j=1,\dots, n{+}1,$ in a $(n+1)$-dimensional
thread module $V=\oplus_{j=1}^{n+1} V_j$, such that
$f_j \in V_j$,
gives us a representation of $\mathfrak{g}$ by lower
triangular matrices. Taking the dual module $V^*$
$$
\left(\rho^*(g)\tilde f\right)(v)=\tilde f(\rho(g)v)
$$
with a basis
$f^j \in V_j^*, \;j=1,\dots, n{+}1$
we will get a representation $\rho^*: \mathfrak{g} \to T_n(\mathbb K)$
by upper triagular matrices.

Or one can change the ordering of the basis of $V$ considering a new basis
$v_1'=v_{n+1},\dots,v_{n+1}'=v_1$.

\begin{proposition}
Let a Massey product $\langle \omega_1,\omega_2,{\dots},\omega_n \rangle$
be defined and trivial in $H^2(\mathfrak{g})$ for some $1$-cohomology classes
$\omega_i \in H^1(\mathfrak{g})$ of a Lie algebra $\mathfrak{g}$. Then
$\langle x_1\omega_1,x_2\omega_2,{\dots},x_n\omega_n \rangle$
is also defined and trivial for any choice of non-zero constants
$x_1,x_2,{\dots},x_n$.
\end{proposition}

\section{Massey products in $H^*(L_1)$ and $H^*(\mathfrak{m}_0)$}

We recall that the algebra $H^*(L_1)$ has a trivial multiplication.
Buchstaber conjectured that the algebra
$H^*(L_1)$ is generated with respect to the Massey products
by $H^1(L_1)$, moreover all corresponding Massey
products can be chosen non-trivial.
The first part of Buchstaber's conjecture was proved
by Feigin, Fuchs and Retakh \cite{FeFuRe}.

\begin{theorem}[\cite{FeFuRe}]
\label{FeigFu}
For any $q \ge 2$
$$
g^q_- \in {\mathbb K}\langle g^{q{-}1}_+,\underbrace{e^1,\dots,e^1}_{2q-1}
\rangle,
\quad
g^q_+ \in {\mathbb K}\langle g^{q{-}1}_+,\underbrace{e^1,\dots,e^1}_{3q-1}
\rangle.
$$
\end{theorem}

The second part of the Buchstaber conjecture was not treated
in \cite{FeFuRe}, moreover it follows from the Proposition
\ref{triviality} that the $2q$-fold Massey product
$\langle g^{q{-}1}_+,\underbrace{e^1,\dots,e^1}_{2q-1}\rangle$
from the Theorem \ref{FeigFu} is trivial. Triviality or non-triviality
of the $3q$-fold Massey products from the Theorem \ref{FeigFu}
is even less clear and have not been studied yet.

In 2000 Buchstaber's PhD-student Artel'nykh
considered the second part of the Buchstaber conjecture.
In particular he claimed the following theorem.

\begin{theorem}[\cite{Artel}]
\label{Th_Artel}
There are non-trivial Massey products
$$
g^q_- \in {\mathbb K}\langle \underbrace{e^2,\dots,e^2}_{q-1},
g^{q{-}1}_+,e^1 \rangle, q \ge 2,
\quad
g^{2l{+}1}_+ \in {\mathbb K}\langle \underbrace{e^2,\dots,e^2}_{3l+1},
g^{2l}_+ \rangle, l \ge 1.
$$
\end{theorem}
One can see that the cohomology classes $g^{2l}_+$ were not
represented by means of non-trivial Massey products. On the another hand
Artel'nykh's article contain only a brief sketch of the proof.
Hence the second part of Buchstaber's conjecture is still open.

We have mentioned that the Massey products in $H^*(\mathfrak{g})$
and $H^*({\rm gr}_C\mathfrak{g})$ are closely related. Recall that
${\rm gr}_C L_1 \cong \mathfrak{m}_0$.
We came to the problem of description of Massey products
in the cohomology $H^*(\mathfrak{m}_0)$. The special question
is the description
of equivalency classes of trivial Massey products
$\langle \omega_1,\dots,\omega_n\rangle$ of $1$-cohomology classes
$\omega_1,\dots,\omega_n$. The purpose of this interest is
to consider
Massey products of the form $\langle \omega_1,\dots,\omega_n,\Omega\rangle$,
where $\Omega$ is an element from $H^*(\mathfrak{g})$.

An infinite dimensional space $H^2(\mathfrak{m}_0)$ is spanned
by the cohomology classes of following $2$-cocycles \cite{FialMill}
\begin{multline}
\label{2-hom}
\omega(e^2{\wedge} e^3)=e^2{\wedge} e^3,
\omega(e^3{\wedge} e^4)=e^3{\wedge} e^4 -e^2{\wedge} e^5,\\
\omega(e^4{\wedge} e^5)=e^4{\wedge} e^5 -e^3{\wedge} e^6+
e^2{\wedge} e^7, \dots,\\
\omega(e^k{\wedge} e^{k{+}1})=
\sum_{l=0}^{k-2} ({-}1)^l e^{k-l}{\wedge} e^{k{+}1{+}l}, \dots
\end{multline}

All of the cocycles (\ref{2-hom}) can be represented as Massey products.
Namely let consider the following matrix of a formal connection
$$
A=\left(\begin{array}{ccccccc}
  0      & e^2 & -e^3 & \dots  & (-1)^ke^k & (-1)^{k{+}1}e^{k{+}1}   & 0  \\
  0      & 0      & e^1 & 0& \dots & 0   & e^{k{+}1}   \\
  0      & 0      & 0 & e^1 &  0 & \dots   & e^k   \\
  \dots & \dots &   \dots     & \dots  &\dots  & \dots & \dots  \\
  0  &0    &  0&   \dots  &  0     &  e^1 & e^3 \\
  0   &0   & 0      & 0      & \dots & 0          & e^2   \\
  0    &0  & 0      & 0      &  0& \dots          & 0
  \end{array}\right).
$$
For the related cocycle $c(A)$ we have
$$
c(A)=\sum_{l=2}^{k{+}1}(-1)^l e^l \wedge e^{k+3-l}=
2\omega(e^k{\wedge}e^{k{+}1}).
$$
Thus we have proved the following
\begin{proposition}
$$
2\omega(e^k{\wedge}e^{k{+}1}) \in
\langle e^2, \underbrace{e^1, \dots, e^1}_{2k-3}, e^2 \rangle,
\quad k \ge 2.
$$
\end{proposition}

\begin{example}
We take $k=2$ and the matrix $A$ of a formal connection
that corresponds to $\langle e^2,e^1,e^2\rangle$
$$
A=\left(\begin{array}{cccc}
0 & e^2 & -e^3& 0\\
0 & 0 & e^1& e^3\\
0 & 0 & 0& e^2\\
0 & 0 & 0& 0 \\
\end{array}\right).
$$
The related cocycle $c(A)$ is equal to
$$
c(A)=2e^2{\wedge} e^3=2\omega(e^2{\wedge} e^3).
$$
\end{example}

The space $H^1(\mathfrak{m}_0)$ is spanned by $e^1$ and $e^2$
and hence every $n$-fold Massey product of
elements from $H^1(\mathfrak{m}_0)$ has a form
$$
\underbrace{\langle \alpha_1 e^1 +\beta_1 e^2, \alpha_2 e^1 +\beta_2 e^2,
\dots, \alpha_n e^1 +\beta_n e^2 \rangle}_n.
$$
A product $e^1 \wedge e^2=de^3$ is cohomologicaly trivial,
hence a triple product
$$\langle \omega_1, \omega_2, \omega_3 \rangle=
\langle \alpha_1 e^1 +\beta_1 e^2, \alpha_2 e^1 +\beta_2 e^2,
\alpha_3 e^1 +\beta_3 e^2 \rangle$$
is defined for any choice
of constants $\alpha_i, \beta_i \in {\mathbb K}, i=1,2,3.$
$$
A=\left(\begin{array}{cccccc}
0 & \omega_1 & \gamma_1e^3& 0\\
0 & 0 & \omega_2& \gamma_2e^3\\
0 & 0 & 0& \omega_3\\
0 & 0 & 0& 0 \\
\end{array}\right), \; \gamma_1=\alpha_1\beta_2-\alpha_2\beta_1,
\gamma_2=\alpha_2\beta_3-\alpha_3\beta_2,
$$
The related cocycle $c(A)=\gamma_2\omega_1{\wedge}e^3-
\gamma_1\omega_3{\wedge}e^3$ is trivial if and only if
\begin{equation}
\label{triple_M}
\beta_1(\alpha_2\beta_3-\alpha_3\beta_2)-
\beta_3(\alpha_1\beta_2-\alpha_2\beta_1)=0.
\end{equation}
We have mentioned above that if $n$-fold Massey product
$\langle \omega_1, \omega_2, {\dots},\omega_n \rangle$
is defined than all $(p+1)$-fold Massey products
$\langle \omega_i, \omega_{i+1}, \dots,\omega_{i+p} \rangle$  for
$1\le i \le n-1, 1\le p \le n-2, i+p \le n$
are defined and trivial.

\begin{theorem}
\label{1_Massey}
Up to an equivalence
the following trivial $n$-fold Massey products of non-zero cohomology
classes from $H^1(\mathfrak{m}_0)$ are defined:

$$
\begin{tabular}{|c|c|c|}
\hline
&&\\[-10pt]
name &  Massey product & parameters\\
&&\\[-10pt]
\hline
&&\\[-10pt]
$A^{n+1}_{\lambda}$
 & $\underbrace{\langle \alpha e^1{+}\beta e^2,\alpha e^1{+}\beta e^2,
\dots,\alpha e^1{+}\beta e^2\rangle}_n$ &
$n \ge 3, \; \lambda=(\alpha,\beta) \in {\mathbb K}P^1$ \\
&&\\[-10pt]
\hline
&&\\[-10pt]
$B^{n+1}_{\alpha,\beta}$&
 $\underbrace{\langle \lambda_1 e^1{+}e^2,\lambda_2 e^1{+}e^2,
{\dots},\lambda_n e^1{+}e^2 \rangle}_n$ &
$\begin{array}{c}
n \ge 3, \; \lambda_i=i \alpha{+}\beta, \\
\alpha,\beta \in {\mathbb K}, \; \alpha \ne 0 \end{array}$
\\
&&\\[-10pt]
\hline
&&\\[-10pt]
$C^{n+1}_{l,\alpha}$ &
$\underbrace{\langle e^1, {\dots}, e^1}_l, e^2{+}\alpha e^1 ,
\underbrace{e^1,\dots, e^1\rangle}_{n-l-1}$ &
$\begin{array}{c} \alpha \in {\mathbb K},\; n \ge 3,\\
0\le l \le n{-}1 \end{array}$  \\
&&\\[-10pt]
\hline
&&\\[-10pt]
$D^{2k+3}_{\alpha,\beta}$&
 $\langle e^2+\alpha e^1, \underbrace{e^1,\dots,e^1}_{2k},
e^2+\beta e^1\rangle$
 & $k \ge 1, \;\alpha,\beta \in {\mathbb K}$ \\[2pt]
\hline
\end{tabular}
$$
\end{theorem}

\begin{proof}
The statement of the present theorem is equivalent to
Benoist's classification \cite{Ben} of indecomposable
thread modules over $\mathfrak{m}_0$
($\mathfrak{m}_0$ is graded as
${\rm gr}_C\mathfrak{m}_0 \cong \mathfrak{m}_0$).
More pricesly we consider a finite-dimensional $\mathfrak{m}_0$-module $V$
with a basis $v_1,\dots,v_{n{+}1}$ such that
$$
e_1 v_i =\alpha_i v_{i-1}, \; e_2 v_i =\beta_i v_{i-1}, \;
i=1,\dots,n+1.
$$
in the last formula we assume that $v_0=0$.
It is sufficient to define only an action of $e_1$ and $e_2$ on $V$ because
the algebra $\mathfrak{m}_0$ is generated by them.

Taking the corresponding matrix $A$ of this representation
with respect to the basis $v_1,\dots,v_{n+1}$ we see that it has
elements $\omega_i=\alpha_i e^1 +\beta_i e^2$ at its second diagonal.
One can regard $A$ as a defining system of the Massey product
$\langle \underbrace{\omega_1,
\dots,\omega_n\rangle}_n, \omega_i=\alpha_i e^1 {+}\beta_i e^2$.
Obviously $V$ is decomposable as a direct sum of thread modules
iff $e_1v_i=e_2v_i=0$ for some $i, 1\le i \le n{+}1$. The last
condition means that $\omega_i=0$ at the second diagonal
of the matrix $A$.

We will prove this theorem by induction and start with triple
products.
First of all let consider the case when all $\beta_i \ne 0, i=1,2,3$.
Taking an equivalent defining system one can asume that $\beta_i=1,\:
i=1,2,3$. In terms of representations it means that one
can choose a base $v_1,\dots,v_4$  of $V$ such as
$$
e_2v_i=v_{i-1}, \quad i=1,\dots, 4.
$$
The equation (\ref{triple_M}) in our case looks in a following way
$$
2\alpha_2-\alpha_1-\alpha_3=0,
$$
and it means that the numbers $\alpha_1,\alpha_2,\alpha_3$
form an arithmetic progression and we have either
a type $B^{4}_{\alpha,\beta}$ or $A^{4}_{\lambda}$ with
$\lambda=(\alpha,1)$. We keep Benoist's notations of types
of thread modules from \cite{Ben}.

The equation (\ref{triple_M}) implies that if one constant from
$\beta_1,\beta_2,\beta_3$ is equal to zero, then at least another
one $\beta_i$ is trivial also.

The case when $\beta_1=\beta_2=\beta_3=0$ is
equivalent to $\langle e^1, e^1, e^1 \rangle$, i.e. to the type $A^4_{1,0}$
from the table above. The remaining three possibilities are
$C^4_{l,\alpha}, \:l=0,1,2$.

Let a $n$-fold product
$\langle \omega_1, \dots, \omega_n \rangle$ be defined and trivial then
$(n-1)$-fold product $\langle \omega_1, \dots, \omega_{n-1} \rangle$
is also trivial and by our inductive assumption is equivalent to some case
from the table above.

The triple Massey product
$\langle \omega_{n-2}, \omega_{n-1}, \omega_{n}\rangle$ is trivial in its turn
and one can write out the equation (\ref{triple_M}) for
all classes $\langle \omega_1, \dots, \omega_{n-1} \rangle$:
$$
\begin{tabular}{|c|c|c|}
\hline
&&\\[-10pt]
$\langle \omega_1, \dots, \omega_{n-1} \rangle$ &
equation (\ref{triple_M}) for
$\langle \omega_{n-2}, \omega_{n-1},\omega_n \rangle$&
$\langle \omega_1, \dots, \omega_n \rangle$\\
&&\\[-10pt]
\hline
&&\\[-10pt]
$A^{n}_{\lambda},\: \lambda=(\alpha,\beta)$
 & $\alpha \beta \beta_n=\beta^2\alpha_n$ &
$A^{n+1}_{\lambda}, C^{n+1}_{n{-}1,\rho}$ \\
&&\\[-10pt]
\hline
&&\\[-10pt]
$B^{n}_{\alpha,\beta}$&
 $(\alpha n+\beta)\beta_n=\alpha_n$ &
$B^{n+1}_{\alpha,\beta}$\\
&&\\[-10pt]
\hline
&&\\[-10pt]
$C^{n}_{l,\alpha}$ &
$\begin{array}{c}
0=0,\; l < n{-}2;\\
2\beta_n=0, \; l=n{-}2;\\
-\beta_n=0,\; l=n{-}1.
\end{array}$ &
$\begin{array}{c}
?\\
C^{n+1}_{l,\alpha}\\
C^{n+1}_{l,\alpha}
\end{array}$  \\
&&\\[-10pt]
\hline
&&\\[-10pt]
$D^{2k+3}_{\alpha,\beta}$&
 $-\beta_n=0$
& ?\\[2pt]
\hline
\end{tabular}
$$
\begin{lemma}
The following Massey products
$$
\langle e^2+\alpha e^1, \underbrace{e^1,\dots,e^1}_{2k},
e^2+\beta e^1, e^1\rangle, \quad
\langle e^1,e^2+\alpha e^1, \underbrace{e^1,\dots,e^1}_{2k},
e^2+\beta e^1\rangle
$$
are defined and non-trivial.
\end{lemma}
\begin{proof}
In order to symplify the notations
we are going to consider only
the first Massey product with $k=1$, the general case can be treated
analogously. We will write only non-trivial
entries of the defining system $A$ (so-called the Massey triangle)
$$
\begin{array}{ccccc}
e^2{+}\alpha e^1 & -e^3{+}\dots  & e^4{+}\dots &(\alpha{-}\beta)e^5{+}\dots & \\
                 & e^1           & 0+\dots     & \dots                      & -3e^5{+}\dots\\
                 &               &  e^1        &  e^3{+}\dots               &  -2e^4{+}\dots\\
                 &               &             &  e^2{+}\beta e^1           & -e^3{+}\dots\\
                 &               &             &                            & e^1 \\
\end{array}
$$
The related cocycle $c(A)$ is equal to
$$
c(A)=3e^3{\wedge} e^4-3e^2{\wedge} e^5+\dots=3\omega(e^3{\wedge} e^4)+\dots,
$$
where dots stand everywhere instead of summands of lower gradings.
\end{proof}
Hence one can coclude that a trivial Massey product of type $D$ is not
extendable to a higher fold trivial Massey product.
The same arguments show that $C^n_{l,\alpha}$ is extended
only to $C^{n{+}1}_{l,\alpha}$ if $ l > 0$ and
we have two possibilities in the case $C^{2k{+}1}_{0,\beta}$.
They are $C^{2k{+}2}_{0,\beta}$ and $D^{k}_{\alpha,\beta}$
respectively.
\end{proof}

\begin{theorem}
\label{main_th}
The cohomology algebra  $H^*(\mathfrak{m}_0)$  is generated
with respect to the non-trivial Massey products
by $H^1(\mathfrak{m}_0)$, namely
\begin{equation}
\begin{split}
\omega(e^{2}{\wedge}e^{i_2}{\wedge}{\dots}
{\wedge} e^{i_q}{\wedge} e^{i_q{+}1})=
e^{2}{\wedge}\omega (e^{i_2}{\wedge}{\dots}
{\wedge} e^{i_q}{\wedge} e^{i_q{+}1}),\\
2\omega(e^{k}{\wedge}e^{k{+}1}) \in
\langle e^2, \underbrace{e^1, \dots, e^1}_{2k-3},e^2 \rangle,\\
(-1)^{i_1}\omega(e^{i_1}{\wedge} e^{i_2}{\wedge} {\dots}
{\wedge} e^{i_q}{\wedge} e^{i_q{+}1}) \in
\langle e^2, \underbrace{e^1, \dots, e^1}_{i_1-2},
\omega(e^{i_2}{\wedge}{\dots}
{\wedge} e^{i_q}{\wedge} e^{i_q{+}1})
\rangle,
\end{split}
\end{equation}
\end{theorem}

\begin{proof}

First of all we present a graded defining system (a graded
formal connection) $A$ for a Massey product
$\langle e^2, e^1, \dots, e^1,\omega(e^{i_2}{\wedge}{\dots}
{\wedge} e^{i_q}{\wedge} e^{i_q{+}1}) \rangle$.
To symplify the formulas we will write $\omega$ instead of
$\omega(e^{i_2}{\wedge}{\dots}
{\wedge} e^{i_q}{\wedge} e^{i_q{+}1})$.

One can verify that the following matrix $A$
with  non-zero entries at the second diagonal, first line and
first row gives us an answer.

$$
A=\left(\begin{array}{ccccccc}
  0      & e^2 & -e^3 & e^4  & \dots & (-1)^{i_1}e^{i_1}   & 0  \\
  0      & 0      & e^1 & 0& \dots & 0   & D_{-1}^{i_1-2}\omega   \\
  0      & 0      & 0 & e^1 &  \dots & 0   & D_{-1}^{i_1-3}\omega   \\

   &      &        &        & \dots &            &           \\
  0  &0    & 0      & 0      & \dots & e^1 & D_{-1}\omega \\
  0   &0   & 0      & 0      & \dots & 0          & \omega   \\
  0    &0  & 0      & 0      & \dots & 0          & 0
  \end{array}\right).
$$
The proof follows from
$$
d(D_{-1}^k \omega)=e^1{\wedge}D_{-1}^{k{-}1} \omega, \quad
d((-1)^k e^k)=(-1)^{k{-}1}e^{k{-}1}{\wedge}e^1.
$$
The related cocycle $c(A)$ is equal to
$$
c(A)=\sum_{l\ge 2}^{i_1}(-1)^l e^{l} \wedge D_{-1}^{i_1{-}l} \omega
=(-1)^{i_1}\sum_{k\ge 0}^{i_1{-}2}(-1)^k D_1^k e^{i_1} \wedge D_{-1}^k \omega
$$

\begin{example}
We take $\langle e^2, e^1, \omega(e^4{\wedge} e^5)\rangle $.
$$
A=\left(\begin{array}{cccc}
  0      & e^2 & -e^3 &  0  \\
  0      & 0      & e^1 & D_{-1}\omega(e^4{\wedge} e^5)   \\
  0   &0    & 0          & \omega(e^4{\wedge} e^5)   \\
  0    &0  & 0      & 0
    \end{array}\right).
$$
One can verify the computations of its related cocycle $c(A)$
\begin{multline}
c(A)=e^2\wedge D_{-1}\omega(e^4{\wedge} e^5)
-e^3 \wedge \omega(e^4{\wedge} e^5)=\\
=e^2\wedge (e^4{\wedge} e^6-2e^3{\wedge} e^7+3e^2{\wedge} e^8)
-e^3 \wedge (e^4{\wedge} e^5-e^3{\wedge} e^6+e^2{\wedge} e^7)=\\
=-e^3 {\wedge} e^4{\wedge} e^5+e^2 {\wedge} e^4{\wedge} e^6
-e^2 {\wedge} e^3{\wedge} e^7=-\omega(e^3 {\wedge} e^4{\wedge} e^5).
\end{multline}
\end{example}

\begin{lemma}
$$
\sum_{k\ge 0}(-1)^k D_1^k e^{i_1} \wedge D_{-1}^k
\omega= \omega(e^{i_1}{\wedge} e^{i_2}{\wedge} {\dots}
{\wedge} e^{i_q}{\wedge} e^{i_q{+}1})
$$
\end{lemma}
\begin{proof}
This lemma is almost evident. Both parts of the equality
are closed forms. One has to compare the monomials
of the form
$e^{j_1}{\wedge} e^{j_2}{\wedge} {\dots}
{\wedge} e^{j_q}{\wedge} e^{j_q{+}1}$ in the decompositions
of the left-hand and right-hand sides.
The operator $D_{-1}$ strictly increases the difference between
two last superscripts of monomials
$$
D_{-1}(e^{j_1}{\wedge} e^{j_2}{\wedge} {\dots}
{\wedge} e^{j_q}{\wedge} e^{j_q{+}1})
= \sum_{k\ge 0}
(-1)^k D_{1}^k (e^{j_1}{\wedge} e^{j_2}{\wedge} {\dots}{\wedge} e^{j_q})
\wedge D_{-1}^{k{+}1} e^{j_q{+}1}.
$$
Hence there is the only one monomial of the form we are looking for
on the left-hand side and the same one on the right-hand side and it is
$e^{i_1}{\wedge} e^{i_2}{\wedge} {\dots}
{\wedge} e^{i_q}{\wedge} e^{i_q{+}1}$.
\end{proof}
\begin{lemma}
\label{non_trivaiality}
Let $\tilde A$ be an arbitrary defining system (the matrix of
a formal connection) for
$\langle e^2, e^1, \dots, e^1,\omega(e^{i_2}{\wedge}{\dots}
{\wedge} e^{i_q}{\wedge} e^{i_q{+}1}) \rangle$.
Then its related cocycle $c(\tilde A)$ is equal to
$$
(-1)^{i_1}\omega(e^{i_1}{\wedge} {\dots}
{\wedge} e^{i_q}{\wedge} e^{i_q{+}1})+\sum_{j_1{<}i_1}
\lambda_{j_1,\dots,j_q}
\omega(e^{j_1}{\wedge} {\dots}{\wedge} e^{j_q}{\wedge} e^{j_q{+}1})
{+}e^1{\wedge} \Omega,
$$
\end{lemma}
for some constants $\lambda_{j_1,\dots,j_q}$ and $q$-form $\Omega$.

\begin{proof}

We will rewrite our defining system $\tilde A$ in the form of a Massey
triangle of the defining system $\tilde A$.
$$
\begin{array}{cccccc}
  e^2 & -e^3{+}\rho_1^1 & e^4{+}\rho_1^2  & \dots & (-1)^{i_1}e^{i_1}{+}\rho_1^{i_1{-}2}   &   \\
      & e^1           & \rho_2^1 & \dots &\dots& D_{-1}^{i_1{-}2}\omega{+}{\dots}{+}\Omega_{i_1-2}   \\
      &               &   \dots  & \dots &\dots&   \dots      \\
      &               &          &  e^1  &\rho_{i_1{-}2}^1&D_{-1}^2\omega{+}{\dots}{+}\Omega_2 \\
      &               &          &       & e^1 & D_{-1}\omega{+}\Omega_1     \\
      &               &          &       &     & \omega                      \\
  \end{array}
$$
$\Omega_i$ are some closed $q$-forms.
$1$-forms
$\rho_1^1,\dots,\rho_{i_1{-}2}^1$ standing at the second diagonal
are also closed and therefore they are linear combinations of $e^1$ and $e^2$.
Continue this procedure and using the Maurer-Cartan equation
and an inductive assumption
it is easy to see that $\rho_k^l$ is a linear combination
of $e^1,\dots,e^{l{+}1}$. Hence the maximal value of superscript that
we can meet at the first line of our triangle is
$i_1$. Thus one can conclude that we have the only one monomial
of the form $e^{i_1}{\wedge}{\dots}{\wedge}e^j{\wedge}e^{j{+}1}$,
$i_1{<}\dots{<}j$
in the decompostion of the related cocycle $c(\tilde A)$. It comes
from the summand $(-1)^{i_1}(e^{i_1}{+}\rho_1^{i_1{-}2}){\wedge}\omega$
in the formula of $c(\tilde A)$ and  it is
$(-1)^{i_1}e^{i_1}{\wedge} {\dots}{\wedge} e^{i_q}{\wedge} e^{i_q{+}1}$.
\end{proof}
The Lemma \ref{non_trivaiality} provides us
with a proof of non-triviality of the Massey products
$\langle e^2, e^1, \dots, e^1,\omega(e^{i_2}{\wedge}{\dots}
{\wedge} e^{i_q}{\wedge} e^{i_q{+}1}) \rangle$.
The proof in the case $\langle e^2,e^1,\dots,e^1,e^2 \rangle$ can be
obtained by the same arguments.
\end{proof}

\end{document}